\documentclass[11pt]{amsart}

\usepackage{amsmath, amsfonts, amssymb}

\newtheorem{thm}{Theorem}[section]
\newtheorem{lem}[thm]{Lemma}
\newtheorem{prop}[thm]{Proposition}
\newtheorem{cor}[thm]{Corollary}

\usepackage{tikz}

\usetikzlibrary{decorations.markings}
\usetikzlibrary{shapes.geometric}

\pgfdeclarelayer{edgelayer}
\pgfdeclarelayer{nodelayer}
\pgfsetlayers{edgelayer,nodelayer,main}

\tikzstyle{none}=[inner sep=0pt]
\definecolor{hexcolor0xf81e1c}{rgb}{0.973,0.118,0.110}
\definecolor{hexcolor0x3c00ff}{rgb}{0.235,0.000,1.000}

\tikzstyle{blackvertex}=[circle,fill=black,draw=black, scale=0.61]
\tikzstyle{vertex}=[circle, fill=white,draw=black, scale=0.61]
\tikzstyle{box}=[rectangle,fill=none,draw=none]

\tikzstyle{arc}=[-,draw=black,postaction={decorate},decoration={markings,mark=at position .8 with {\arrow{>}}},line width=0.61]

\begin{document}
\title[Complexity of inj.~hom.~to small tournaments, and to inj.~oriented colourings]{Complexity of injective homomorphisms to small tournaments, and of injective oriented colourings}

\author[R.J.~Campbell]{Russell J.~Campbell} \address{Russell J.~Campbell, University of the Fraser Valley\\ Abbotsford, British Columbia}
\author[N.E.~Clarke]{Nancy E.~Clarke$^\ast$} \address{Nancy E.~Clarke, Acadia University \\ Wolfville, Nova Scotia} \thanks{$^\ast$ Research supported by NSERC}
\author[G.~MacGillivray]{Gary MacGillivray$^\ast$}\address{G.~MacGillivray, University of Victoria\\ Victoria, British Columbia} 


\begin{abstract}
Several possible definitions of local injectivity for a homomorphism of an oriented graph $G$ to an oriented graph $H$ are considered.  
In each case, we determine the complexity of deciding whether there exists such a homomorphism when $G$ is given and $H$ is a fixed  tournament on three or fewer vertices.
Each possible definition leads to a locally-injective oriented colouring problem.
A dichotomy theorem is proved in each case. 
\end{abstract}

\maketitle

\small
\begin{center}
\textbf{Keywords:} injective graph homomorphism, oriented colouring, complexity

\textbf{Mathematics Subject Classifications (2010):} 05C15, 05C60, 05C85
\end{center}

\normalsize

\section{Introduction}
\label{IntroSec}

Three natural possible definitions of local injectivity of a homomorphism $f$ from an \emph{input} oriented graph $G$ to a \emph{target} oriented graph $H$ are: for every vertex $x \in V(G)$, the function $f$ is injective when restricted to:
\begin{enumerate}
\item the  in-neighbourhood $N^-(x) $; or
\label{PropIO}
\item $N^-(x)$ and $N^+(x) $ separately; or
\label{PropIOS}
\item the union $N^-(x) \cup N^+(x)$.
\label{PropIOT}
\end{enumerate}
When $H$ is \emph{reflexive}, that is,  has a loop at every vertex, the three definitions are different.  When $H$ is \emph{irreflexive}, that is,  has no loops, definitions \ref{PropIOS} and \ref{PropIOT} coincide.   Each of these five situations leads naturally to a notion of locally-injective oriented $k$-colouring.

Locally-injective homomorphisms (as in possible definition \ref{PropIO}) and colourings of oriented graphs were first introduced as an example in monadic second order logic \cite{courcelle}.  Consequently, by Courcelle's Theorem, these problems are all solvable in polynomial time when the input has bounded treewidth.  The same holds for the other possible definitions above.  

Possible definition \ref{PropIO} has been studied in previous papers for both irreflexive and reflexive targets \cite{mrs,mrs2,mrs3,mrs1,cobusthesis}.  A fairly complete theory has been developed.  When the target, $H$,  is reflexive there is a dichotomy theorem characterizing the oriented graphs $H$ for which the problem of deciding the existence of a homomorphism to $H$ is Polynomial, and those for which it is NP-complete.  When $H$ is irreflexive the complexity has been determined when $H$ has maximum in-degree $\Delta^- \geq 3$ or $\Delta^- \leq 1$; when $\Delta^-=2$ the situation is as rich as that for all digraph homomorphism problems, and hence all constraint satisfaction problems \cite{mrs1}.  

Possible definitions \ref{PropIOS} and \ref{PropIOT} have been studied in \cite{russell, CCM}.  Obstructions to
(subgraphs that prevent the existence of) homomorphisms to small tournaments are the focus of \cite{CCM}.  Both definitions 
are considered.
Possible definition \ref{PropIOT} is the main focus of \cite{russell}.

Locally-injective colourings of undirected graphs were first explicitly studied by Hahn, Kratochvil, Si\v{r}an and Sotteau \cite{gena}.  Subsequent papers have considered chordal graphs \cite{hrs08}, planar graphs (see \cite{KS}) and other graph classes, as well as list versions \cite{fk06}.  The complexity of locally-injective homomorphisms has been extensively studied by Fiala, Kratochvil, and others (e.g. see \cite{fkp08,fpt08}).  

The purpose of this paper is to contribute to the theory of locally-injective homomorphisms and colourings under possible definitions \ref{PropIOS} and \ref{PropIOT} above.   In each of the three cases that arise, the complexity of deciding the existence of a homomorphism to $H$ is determined for the four tournaments on at most three vertices.  These results appear in Sections 3, 4, and 5.   Later, in Section 6, these results are then used to determine the complexity of the associated locally-injective oriented colouring problems.  

We conclude this section by noting that the complexity of deciding whether a given directed graph $G$ has a homomorphism to a tournament $H$ has been studied \cite{BHM}.
There is a dichotomy theorem: the problem is Polynomial when $H$ has at most one directed cycle, and NP-complete when $H$ has at least two directed cycles.
The results reported in this paper are first steps towards finding a similar theorem for locally-injective homomorphisms.



\section{Notation and terminology}
\label{Sec2}

An \emph{oriented graph} is a directed graph $G$ with the property that for any two different vertices $x$ and $y$, at most one of the arcs $xy, yx$ belongs to $E(G)$. 
An oriented graph $G$ can be viewed as arising from a simple graph $H$ by assigning a direction, or \emph{orientation}, to each edge.
The graph $H$ is called the \emph{underlying graph} of $G$, and $G$ is referred to as \emph{an orientation} of $H$.
The \emph{converse} of an oriented graph $G$ is the oriented graph $G^c$ with the same vertex set as $G$, and arc set $\{yx: xy\in E(G)\}$.  

An oriented graph is \emph{reflexive} if it has a loop at each vertex, and \emph{irreflexive} if it has no loops.  The superscript ``$r$'', as in $C_3^r$, indicates that the oriented graph under consideration is reflexive.  
Oriented graphs without this superscript, as in $G$, are irreflexive.

We use $P_n, C_n$, and $T_n$  to denote the directed path on $n$ vertices, the directed cycle on $n$ vertices, and the transitive tournament on $n$ vertices, respectively, $n \geq 1$.  It will be assumed throughout that $C_3$ has vertex set $\{c_1, c_2, c_3\}$ and arc set $\{c_1c_2, c_2c_3, c_3c_1\}$, and that $T_n$ has vertex set $\{t_0, t_1, \ldots, $ $t_{n-1}\}$ and arc set $\{t_it_j: i < j\}$. 

A \emph{homomorphism} of an oriented graph $G$ to an oriented graph $H$ is a function $f: V(G) \to V(H)$ such that $f(x)f(y) \in E(H)$ whenever $xy \in E(G)$.  
When $H$ has a loop, any directed graph has a homomorphism to $H$: map all vertices of $G$ to a vertex of $H$ with a loop.  
Thus, when loops are present, the existence of a homomorphism is a non-trivial question only in the presence of some side condition like selecting the image of each vertex from a list of possible images, or local injectivity.  
The book \cite{hn_book} contains a wealth of information about homomorphism of graphs and digraphs.

We call a homomorphism $f$ of an oriented graph $G$ to an oriented graph $H$:
\begin{itemize}
\item \emph{ios-injective} if, for every vertex $x$ of $G$, the restriction of $f$ to $N^-(x)$ is injective, as is the restriction of $f$ to $N^+(x)$; and
\item   \emph{iot-injective} if, for every vertex $x$ of $G$, the restriction of $f$ to $N^-(x) \cup N^+(x)$ is injective.
\end{itemize}
These two concepts are the same when $H$ is an irreflexive oriented graph, and different when $H$ is a reflexive oriented graph.

The designations ``ios'' and ``iot'' arise from the local injectivity being on \underline{{\bf i}}n-neighbourhoods and \underline{{\bf o}}ut-neighbourhoods \underline{{\bf s}}eparately, and on \underline{{\bf i}}n-neighbourhoods and \underline{{\bf o}}ut-neighbourhoods \underline{{\bf t}}ogether.  In introducing the designations ``ios'' and ``iot'',  the qualifier ``locally'' has been dropped as it is part of the definition. 


It is easy to see that the composition of two ios-injective homomorphisms is an ios-injective homomorphism, and similarly for iot-injective homomorphisms. 


The following structure and its converse will be particularly useful.
We define the \textit{hat} $H_3$ to be the oriented graph with vertex set $V(H_3) = \{v_0, v_1, v_2\}$ 
and edge set $E(H_3) = \{v_0v_1, v_2v_1\}$.   The vertices $v_0$ and $v_2$ will be referred to as the \emph{ends} of $H_3$ or $H_3^c$.  Whether or not $H$ is reflexive, in an ios-injective or iot-injective homomorphism of $H_3$ or $H_3^c$ to $H$, the vertices $v_0$ and $v_2$ must have different images.

\begin{figure}[htbp]
\begin{center}
\begin{tikzpicture}
	\begin{pgfonlayer}{nodelayer}
		\node [style=vertex] (0) at (-5, 1) {};
		\node [style=vertex] (1) at (-4.5, 2) {};
		\node [style=vertex] (2) at (-4, 1) {};
		\node [style=vertex] (3) at (-2, 1) {};
		\node [style=vertex] (4) at (-1.5, 2) {};
		\node [style=vertex] (5) at (-1, 1) {};
		\node [style=vertex] (6) at (-4.5, 2) {};
		\node [style=box] (7) at (-4.5, 0.25) {$H_3$};
		\node [style=box] (8) at (-1.5, 0.25) {$H_3^c$};
	\end{pgfonlayer}
	\begin{pgfonlayer}{edgelayer}
		\draw [style=arc] (0) to (1);
		\draw [style=arc] (2) to (1);
		\draw [style=arc] (4) to (3);
		\draw [style=arc] (4) to (5);
	\end{pgfonlayer}
\end{tikzpicture}
\caption{The hat and its converse.}
\label{fig1}
\end{center}
\end{figure}
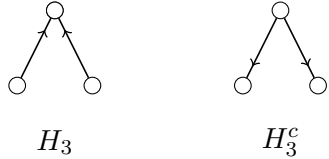

\section{Irreflexive targets}
In this section we show that, if $T$ is an irreflexive tournament on at most 3 vertices, then the problem of deciding whether a
given oriented graph has an ios-injective (and hence also iot-injective) homomorphism to $T$ is Polynomial.
A given oriented graph has an ios-injective homomorphism to $T_1$ if and only if it has no edges, 
and has an ios-injective homomorphism to $T_2$ if and only if it is a disjoint union of copies of $T_1$ and $T_2$.
A given oriented graph, $G$, has an ios-injective homomorphism to $C_3$ if and only if it has maximum in-degree 1, maximum out-degree 1, and has a homomorphism to $C_3$. It follows that $G$ has an ios-injective homomorphism to $C_3$ if and only if it is a disjoint union of directed
paths, and directed cycles of length a multiple of 3. These conditions are easy to check in polynomial time.  It remains to consider ios-injective homomorphisms to the transitive triple.

\begin{prop}
The problem of deciding whether a given oriented graph has an ios-injective homomorphism to $T_3$ is Polynomial.
\end{prop}

\noindent \textbf{Proof.}
Let $G$ be a given digraph.  
If the underlying graph of $G$ has a vertex of degree 3 or more, then $G$ has no ios-injective homomorphism to $T_3$.  
Hence assume that $G$ is an orientation of a graph with maximum degree at most 2.
Therefore the underlying graph of $G$ is a disjoint union of paths and cycles, and hence has treewidth at most 2.
Since ios-injective homomorphism is expressible in monadic second-order logic, the statement now follows from 
Courcelle's Theorem.
 \hfill $\square$

\section{ios-injective homomorphisms to small reflexive targets}

In this section, we determine the complexity of deciding whether there exists an  ios-injective homomorphism from a given oriented graph $G$ to the fixed oriented graph $H$ when $H$ is one of the four reflexive tournaments on at most three vertices.  

It is clear that an oriented graph $G$ has an ios-injective homomorphism to $T_1^r$ if and only if it has maximum in-degree at most one and maximum out-degree at most one, that is, if and only if neither $H_3$ nor $H_3^c$ is a subgraph of $G$.  Consequently, the only oriented graphs which have an ios-injective homomorphism to $T_1^r$ are disjoint unions of directed paths and directed cycles.

\begin{prop}
The problem of deciding whether a given oriented graph has an ios-injective homomorphism to $T_2^r$ is Polynomial.
\label{Prop:iosT_2^r}
\end{prop}

\noindent \textbf{Proof.}
We describe a reduction to 2-SAT.  
Associate the vertices $t_0$ and $t_1$ of $T_2^r$ with false and true, respectively.
Given an oriented graph $G$, the corresponding instance of 2-SAT has the set of variables $\{x_v: v \in V(G)\}$.
Since no oriented graph with a vertex of in-degree at least 3, or a vertex of out-degree at least 3, has an
ios-injective homomorphism to $T_2^r$, we can assume that $\Delta^+(G) \leq 2$ and $\Delta^-(G) \leq 2$.

The set of clauses is constructed as follows.
\begin{itemize}
\item[(i)] If $\mathrm{deg}^+(v) = 2$, then $\neg x_v$ is a clause.
\item[(ii)] If $\mathrm{deg}^-(v) = 2$, then $x_v$ is a clause.
\item[(iii)] If $vw \in E$, then $\neg x_v \vee x_w$ is a clause.
\item[(iv)] If $v$ and $w$ are the ends of a copy of $H_3$ or $H_3^c$, then $ x_v \vee x_w$ and $\neg x_v \vee \neg x_w$ are clauses.
\end{itemize}

All clauses in groups (i) and (ii) are satisfied if and only if the image of any vertex of out-degree 2 is $t_0$ and the image of any vertex of in-degree 2 is $t_1$.
All clauses in group (iii) are satisfied if and only if the mapping corresponding to the truth assignment preserves arcs.
And finally, all clauses in group (iv) are satisfied if and only if the ends of a copy of $H_3$ or $H_3^c$ are assigned different images.
It follows that there is an ios-injective homomorphism of $G$ to $T_2^r$ if and only all clauses are satisfied.
\hfill $\square$

We now show that the problem of deciding the existence of an ios-injective homomorphism to  $C_3^r$ is NP-complete.
Some ``gadget'' oriented graphs which map to $C_3^r$ only in special ways will be used in the NP-completeness proof.
For an integer $d \geq 1$, the oriented graph $D_d$ is constructed from a directed cycle $v_1, v_2, \ldots, v_{6d}, v_1$ by adding the vertices $x_1, x_2, \ldots, x_{3d}$ and arcs $v_{2t}x_t, x_tv_{2t-1},$ $t = 1, 2, \ldots, 3d$.  

\begin{lem}
In an ios-injective homomorphism of $D_d$ to $C_3^r$ the vertices $x_1, x_4,$ $\ldots,$ $ x_{3d-2}$ all have the same image. 
\label{LemmaD_d}
\end{lem}

\noindent \textbf{Proof.}
Let $f$ be an ios-injective homomorphism of $D_d$ to $C_3^r$.  Without loss of generality, suppose $f(v_1) = c_1$.  Then $f(v_2)$ is either $c_1$ or $c_2$.  

Suppose first that $f(v_2) = c_1$.  Then, observing that an ios-injective homomorphism of an irreflexive directed 3-cycle to $C_3^r$ either assigns every vertex the same image, or assigns no two vertices the same image, it must be that $f(x_1) = c_1$.  By injectivity $f(v_3) \not= f(x_1)$, so $f(v_3) = c_2$, the only other out-neighbour of $c_1$.  It follows that $f(x_2) = c_2$.  Similarly, $f(v_4) \neq c_3$, so that $f(v_4) = f(x_2) = c_2$. Continuing in this way, the vertices $v_1, v_2, \ldots, v_{6d}$ map to $c_1, c_1, c_2, c_2, c_3, c_3, c_1, c_1, \ldots, c_3, c_3$, respectively, and the vertices $x_1, x_2, \ldots, x_{3d}$ map to $c_1, c_2, c_3, c_1, \ldots, c_3$, respectively.

Now suppose that $f(v_2) = c_2$.  By our observation regarding homomorphisms of irreflexive directed 3-cycles, it must be that $f(x_1) = c_3$.  Arguing  as in the previous paragraph, ios-injectivity implies $f(v_3) = c_2$,  and $f(x_2) = c_1$, which in turn implies $f(v_4) = c_3$.  Continuing in this way, the vertices $v_1, v_2, \ldots, v_{6d}$ map to  $c_1, c_2, c_2, c_3, c_3, c_1, c_1, \ldots, c_3, c_3, c_1$, respectively, and the vertices $x_1, x_2, \ldots,$ $x_{3d}$ map to $c_3, c_1, c_2, c_3, c_1, \ldots$, $c_2$, respectively. 
\hfill $\square$

\medskip
For $d \geq 2$, let $X_d$ be the oriented graph constructed from $D_d$ by adding $d$ new vertices $n_1, n_2, \ldots, n_d$ and the arcs belonging to $\{x_{3i-2}n_i, n_ix_{3i+1}: i = 1, 2, \ldots, d,\}$, where addition is modulo $3d$.  The following is a consequence of Lemma \ref{LemmaD_d}.

\begin{cor}
In an ios-injective homomorphism of $X_d$ to $C_3^r$, the vertices of the directed cycle $x_1, n_1, x_4,$ $n_2, \ldots, n_d, x_1$ must all be assigned the same image.  Futher, any partial mapping in which these vertices are all assigned the same image can be extended to an ios-injective homomorphism of $X_d$ to $C_3^r$.  
\end{cor}

\begin{thm}
The problem of deciding if a given oriented graph $G$ has an ios-injective homomorphism to  $C_3^r$ is NP-complete.
\label{C_3^rNPc}
\end{thm}

\noindent \textbf{Proof.}
The transformation is  from 3-colouring of graphs with minimum degree at least 3.  Suppose  a graph $G$ is given.  For each vertex $x \in V(G)$,  regard the edges incident with $x$ as being in 1--1 correspondence with the integers $1, 2, \ldots, \mathrm{deg}_G(x)$ so that it is meaningful to talk about the $i^{\textup{th}}$ edge incident with $x$.  Construct a digraph $G^\prime$ as follows.  For each vertex $x \in V(G)$ there is a copy of $X_{\mathrm{deg}_G(x)}$.  (Note that $\mathrm{deg}_G(x) \geq 3$.)  Each edge of $G$ is replaced by an oriented path on three vertices.  Suppose $wz \in E(G)$ is the $i^{\textup{th}}$ edge incident with $w$ and the $j^{\textup{th}}$ edge incident with $z$. Add a new vertex $u_{wz}$ and arcs from  vertex $n_i$ of the copy of $X_{\mathrm{deg}_G(w)}$ corresponding to $w$, and from vertex $n_j$ of the copy of $X_{\mathrm{deg}_G(z)}$ corresponding to $z$, to $u_{wz}$.  The transformation can be accomplished in polynomial time. We will show that $G$ is 3-colourable if and only if there is an ios-injective homomorphism of $G^\prime$ to $C_3^r$.

Suppose that $G$ is 3-colourable, and fix a 3-colouring using the colours $c_1, c_2, c_3$.  If the colour of $x$ is $c_i$, then map vertices $x_1, n_1, x_4, n_2,  \ldots,  n_d, x_1$ of the copy of $X_{\mathrm{deg}_G(x)}$ corresponding to $x$ to $c_i$ and extend this to an ios-injective homomorphism to $C_3^r$.  The ends of each oriented path that replaced an edge of $G$ are now assigned different images, and the mapping so far can be extended to the remaining vertex of each oriented path that replaced an edge of $G$.

Suppose $G^\prime$ has an ios-injective homomorphism to $C_3^r$.  Then, in each copy of $X_{\mathrm{deg}_G(x)}$, all vertices of the directed cycle $x_1, n_1, x_4, n_2, \ldots, n_{\mathrm{deg}_G(w)}, x_1$ are assigned the same image.  Assign this colour to $x$.  By the construction of $G^\prime$ and ios-injectivity, adjacent vertices of $G$ are assigned different colours.
\hfill $\square$

\medskip
We conclude this section by showing that the problem of deciding whether a given oriented graph $G$ has an ios-injective homomorphism to $T_3^r$ is NP-complete.   A useful technical lemma is established first.

\begin{figure}[htbp]
\begin{center}
\begin{tikzpicture}
	\begin{pgfonlayer}{nodelayer}
		\node [style=vertex] (0) at (0, 4.25) [label=above:$d$] {};
		\node [style=vertex] (1) at (3.5, 3.5) [label=above left:$t_2$]{};
		\node [style=vertex] (2) at (4.25, 0.5) [label=right:$c$]{};
		\node [style=vertex] (3) at (3.5, -2.5) [label=above right:$t_0$]{};
		\node [style=vertex] (4) at (0, -3.25) [label=below:$d$] {};
		\node [style=vertex] (5) at (-3.5, -2.5) [label=above left:$t_2$] {};
		\node [style=vertex] (6) at (-4.25, 0.5) [label=left:$c$]{};
		\node [style=vertex] (7) at (-3.5, 3.5) [label=above right:$t_0$] {};
		\node [style=vertex] (8) at (4.25, 4.25) {};
		\node [style=blackvertex] (9) at (4.25, -3.25) [label=right:$v$] {};
		\node [style=vertex] (10) at (-4.25, -3.25) {};
		\node [style=blackvertex] (11) at (-4.25, 4.25) [label=left:$u$] {};
		\node [style=vertex] (12) at (0, 3.25) [label=above right:$t_1$] {};
		\node [style=vertex] (13) at (1.5, 1.25) [label=right:$t_0$]{};
		\node [style=vertex] (14) at (2, 0.5)  {};
		\node [style=vertex] (15) at (1, 0.5) {};
		\node [style=vertex] (16) at (1.5, -0.25) [label=right:$t_2$] {};
		\node [style=vertex] (17) at (1.5, 2.75) [label=above right:$e$] {};
		\node [style=vertex] (18) at (1.5, -1.75) [label=below right:$e$] {};
		\node [style=vertex] (19) at (0, -2.25) [label=below left:$t_1$] {};
		\node [style=vertex] (20) at (-1.5, -1.75) [label=below left:$f$]{};
		\node [style=vertex] (21) at (-1.5, -0.25) [label=left:$t_2$] {};
		\node [style=vertex] (22) at (-1, 0.5) {};
		\node [style=vertex] (23) at (-2, 0.5) {};
		\node [style=vertex] (24) at (-1.5, 1.25) [label=left:$t_0$] {};
		\node [style=vertex] (25) at (-1.5, 2.75) [label=above left:$f$] {};
		\node [style=vertex] (26) at (-3, 0.5) [label=below left:$t_1$] {};
		\node [style=vertex] (27) at (-2.25, 2) [label=above:$a$] {};
		\node [style=vertex] (28) at (-0.75, 2) [label=above:$t_2$] {};
		\node [style=vertex] (29) at (0, 2.5) {};
		\node [style=vertex] (30) at (0, 1.5) {};
		\node [style=vertex] (31) at (0.75, 2) [label=above:$t_0$] {};
		\node [style=vertex] (32) at (2.25, 2)  [label=above:$a$] {};
		\node [style=vertex] (33) at (3, 0.5) [label=below right:$t_1$] {};
		\node [style=vertex] (34) at (2.25, -1) [label=below:$b$] {};
		\node [style=vertex] (35) at (0.75, -1) [label=below:$t_0$] {};
		\node [style=vertex] (36) at (0, -0.5) {};
		\node [style=vertex] (37) at (0, -1.5) {};
		\node [style=vertex] (38) at (-0.75, -1) [label=below:$t_2$] {};
		\node [style=vertex] (39) at (-2.25, -1) [label=below:$b$] {};
	\end{pgfonlayer}[latex arrows]
	\begin{pgfonlayer}{edgelayer}
		\draw [style=arc] (7) to (0);
		\draw [style=arc] (0) to (1);
		\draw [style=arc] (2) to (1);
		\draw [style=arc] (3) to (2);
		\draw [style=arc] (3) to (4);
		\draw [style=arc] (4) to (5);
		\draw [style=arc] (6) to (5);
		\draw [style=arc] (7) to (6);
		\draw [style=arc] (7) to (11);
		\draw [style=arc] (8) to (1);
		\draw [style=arc] (3) to (9);
		\draw [style=arc] (10) to (5);
		\draw [style=arc] (0) to (12);
		\draw [style=arc] (12) to (17);
		\draw [style=arc] (13) to (17);
		\draw [style=arc] (13) to (15);
		\draw [style=arc] (13) to (14);
		\draw [style=arc] (15) to (16);
		\draw [style=arc] (14) to (16);
		\draw [style=arc] (18) to (16);
		\draw [style=arc] (20) to (21);
		\draw [style=arc] (22) to (21);
		\draw [style=arc] (23) to (21);
		\draw [style=arc] (24) to (23);
		\draw [style=arc] (24) to (22);
		\draw [style=arc] (24) to (25);
		\draw [style=arc] (25) to (12);
		\draw [style=arc] (4) to (19);
		\draw [style=arc] (26) to (6);
		\draw [style=arc] (26) to (27);
		\draw [style=arc] (27) to (28);
		\draw [style=arc] (29) to (28);
		\draw [style=arc] (30) to (28);
		\draw [style=arc] (31) to (29);
		\draw [style=arc] (31) to (30);
		\draw [style=arc] (31) to (32);
		\draw [style=arc] (33) to (32);
		\draw [style=arc] (33) to (2);
		\draw [style=arc] (34) to (33);
		\draw [style=arc] (35) to (34);
		\draw [style=arc] (35) to (36);
		\draw [style=arc] (35) to (37);
		\draw [style=arc] (36) to (38);
		\draw [style=arc] (37) to (38);
		\draw [style=arc] (39) to (38);
		\draw [style=arc] (39) to (26);
		\draw [style=arc] (20) to (19);
		\draw [style=arc] (19) to (18);
	\end{pgfonlayer}
\end{tikzpicture}
\caption{The oriented graph $F$ in Lemma \ref{ios-gadget}.}
\label{figF}
\end{center}
\end{figure}
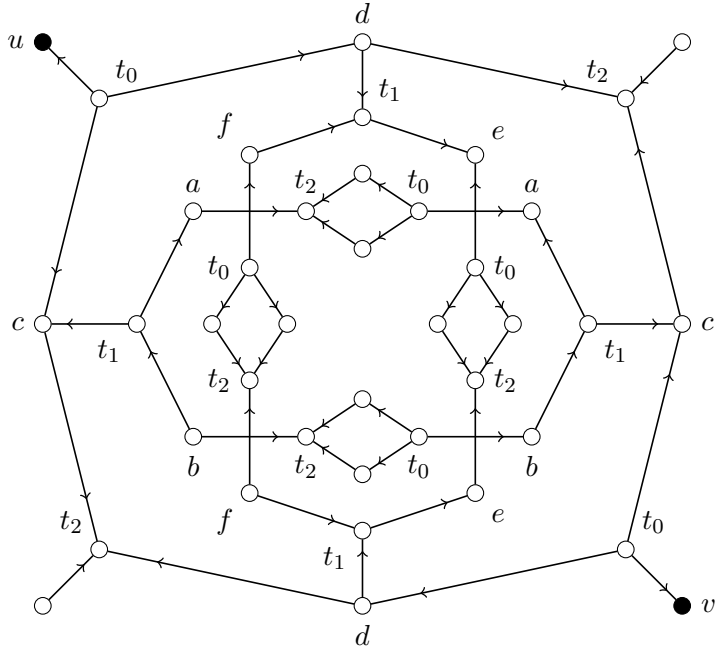


\begin{lem}
Let $F$ be the oriented graph in Figure \ref{figF}.  Then
for $x \in \{t_0, t_1, t_2\}$, there exists an ios-injective homomorphism of $F$ to $T_3^r$ that maps $u$ to $x$, and any such homomorphism also maps $v$ to $x$.
\label{ios-gadget} 
\end{lem} 

\noindent \textbf{Proof.}
We sketch the proof that in an ios-injective homomorphism of $F$ to $T_3^r$ that maps $u$ to $t_1$, the vertex $v$ also maps to $t_1$.  

Referring to Figure \ref{figF}, it is straightforward to check that in any  ios-injective homomorphism of $F$ to $T_3^r$, 
the vertices labelled $t_0,  t_2$ must map to $t_0,  t_2$, respectively. 
It is also easy to check that the vertices labelled $a$ must have the same image,
and similarly for the vertices labelled $b, e$ and $f$.
It will follow from the argument below that the vertices labelled $c$ must have the
same image, and similarly for the vertices labelled $d$.

We show that the vertices labelled $t_1$ must map to $t_1$.
Suppose $u$ maps to $t_1$.  
Then by injectivity its out-neighbour labelled $c$ maps to $t_0$ or $t_2$.
Since $c$ has in-degree 2, it must map to $t_2$.
Therefore $d$ maps to $t_0$.
The in-neighbour of $c$ labelled $t_1$ must map to $t_1$ or $t_2$.
But its in-neighbour labelled $b$ has an out-neighbour labelled $t_2$, so
the in-neighbour of $c$ labelled $t_1$ must map to $t_1$.
By injectivity, the out-neighbour labelled $a$ of this vertex must map to $t_1$,
so the symmetrically located vertex labelled $a$ must also map
to $t_1$, and its in-neighbour labelled $t_1$ must map to $t_1$.
A similar argument shows that the other vertices labelled $t_1$ 
must map to $t_1$.

We now show that $v$ maps to $t_1$.
By the above argument and injectivity, the vertex labelled $c$ on the right of the figure maps to $t_2$.
A symmetric argument shows that the vertex labelled $d$ on the bottom of the figure
must map to $t_0$.
Now, by injectivity, $v$ maps to $t_1$, as wanted.

Similar arguments show that if $u$ maps to $t_0$ then so does $v$, and if $u$ maps to $t_2$ then so does $v$.
\hfill $\square$


\medskip

\begin{thm}
The problem of deciding if a given oriented graph has an ios-injective homomorphism to $T_3^r$ is NP-complete.
\label{iosT_3^r}
\end{thm}

\noindent \textbf{Proof.}
The transformation is from 3-edge colouring of cubic graphs  \cite{holyer}.  Suppose such a graph $G$ is given.  Construct a graph $G^\prime$ as follows.
For each $x \in V(G)$, regard the edges incident with $x$ as being in 1--1 correspondence with the integers $1, 2, 3$ so that it is meaningful to talk about the $i^{\textup{th}}$ edge incident with $x$.  

Let $H_4$ denote the orientation of $K_{1, 3}$ in which there is a vertex of in-degree 3.  In the sequel we refer to $H_4$ as an \emph{in-star}.
Start with a collection of $|V(G)|$ disjoint copies of $H_4$.  Let $S_x$ denote the copy of $H_4$ corresponding to vertex $x$.  Regard the leaves of each oriented graph $S_x$ to be in 1--1 correspondence with $\{1, 2, 3\}$.
Suppose $xy \in E(G)$ is the $i^{\textup{th}}$ edge incident with $x$ and the $j^{\textup{th}}$ edge incident with $y$.  Add a new copy of the oriented graph $F$ shown in Figure \ref{figF} and identify the vertices labelled $u$ and $v$ having in-degree
one with the $i^{\textup{th}}$ leaf of $S_x$ and the $j^{\textup{th}}$ leaf of $S_y$.    The transformation may be accomplished in polynomial time.  We claim that $G$ is $3$-edge-colourable if and only if $G^\prime$ has an ios-injective homomorphism to $T_3^r$.  

Suppose $G$ has a $3$-edge-colouring $f:E(G)\rightarrow \{t_0, t_1, t_2\}$ (the colours are the vertices of $T_3^r$).  For any edge $xy$ of $G$, map the vertices labelled $u$ and $v$ in the corresponding copy of $F$ to $f(xy)$.  Finally, map the centre of each in-star  of $G^\prime$ to its only possible image, $t_2$.  

Conversely, suppose $G^\prime$ has an ios-injective homomorphism to $T_3$.  For each edge $xy$ of $G$, the vertices labelled $u$ and $v$ in the corresponding copy of $F$ in $G^\prime$ must have the same image.  Use this for the colour of $xy$.  The resulting assignment is a 3-edge-colouring because the leaves of each in-star $S_x$ in $G^\prime$ must have different images.
\hfill $\square$

\section{iot-injective homomorphisms to small reflexive targets}

In this section we consider the complexity of deciding whether there exists an iot-injective homomorphism from a given oriented graph $G$ to the fixed oriented graph $H$, when $H$ is one of the four reflexive tournaments on at most three vertices. 

It is clear that an oriented graph has an iot-injective homomorphism to $T_1^r$ if and only if it contains no oriented path on three vertices, that is, if and only if it is a disjoint union of copies of $T_1$ and $T_2$.   
  
We now turn our attention to $T_2^r$.   No orientation of a graph with a vertex of degree three has an iot-injective homomorphism to $T_2^r$.  
Thus, if $G$ admits an iot-injective homomorphism to $T_2^r$, then the underlying graph of $G$ is a disjoint union of paths and cycles.
The following proposition can be proved using a reduction to 2-SAT, or by an appeal to Courcelle's Theorem.

\begin{prop}
The problem of deciding whether a given oriented graph has an iot-injective homomorphism to $T_2^r$ is Polynomial.
\end{prop}

We next consider iot-injective homomorphism to $C_3^r$.
Consider the family of oriented cycles $\mathcal{B}$ such that
each $B \in \mathcal{B}$ is comprised of two disjoint perfect matchings oriented in opposite directions; that is,  $V(B) = \{v_0, v_1, \ldots, v_{2k-1}\}$
and  $E(B) = \{v_0v_1, v_2v_3, \ldots, v_{2k-2}v_{2k-1}\} 
\cup \{v_0v_{2k-1}, v_2v_1,$ $ \ldots, v_{2k-2}v_{2k-3}\}$. 

\begin{lem} 
Let $B \in \mathcal{B}$ have order $n$. Then
(i) $B$ has an iot-injective homomorphism to $C_3^r$ if and only if $n \equiv 0$ (mod 6), and
(ii) $B$ has an iot-injective homomorphism to $T_3^r$ if and only if $n \equiv 0$ (mod 4).
\label{hmm} 
\end{lem}

\noindent{\bf Proof.} 
Let $B \in \mathcal{B}$. 

We first consider iot-injective homomorphism of $B$ to $C_3^r$.
Suppose $B$ has $n$ vertices.
Let $x$ be a vertex of out-degree two.
Without loss of generality $x$ maps to $c_1$.
Then its out-neighbours map to $c_1$ and $c_2$.
Let $y$ be the out-neighbour that maps to $c_1$.
Its out-neighbour must map to $c_2$.
Continuing in this way, starting from $x$, the images of
consecutive vertices are $c_1, c_1, c_2, c_2, c_3, c_3, c_1, c_1, \ldots$.
Therefore an iot-injective homomorphism exists if and only
if f $n \equiv 0$ (mod 6).

We now consider iot-injective homomorphism of $B$ to $T_3^r$.
Suppose $B$ has $n$ vertices.
Let $x$ be a vertex of out-degree two.
Then $x$ maps to $t_0$ or $t_1$.

Suppose first that $x$ maps to $t_0$.
Let $v$ be an out-neighbour of $x$.
If $v$ were mapped to $t_0$, then its other in-neighbour must also map to $t_0$,
in violation of injectivity.
Therefore, the out-neighbours of $x$ map to $t_1$ and $t_2$.
Let $y$ be the out-neighbour that maps to $t_2$.
Then $y$'s other in-neighbour, $z$,  must map to $t_1$
and $z$'s other out-neighbour, $a$, must also map to $t_1$.
The vertex $a$ has another in-neighbour, $b$.
By injectivity, $b$ maps to $t_0$.
Continuing in this way, starting from $x$, the images of consecutive vertices
are $t_0, t_2, t_1, t_1, t_0, \ldots t_0, t_2, t_1, t_1, t_0$.
Therefore $n \equiv 0$ (mod 4).

Now suppose $x$ maps to $t_1$.
As above, the out-neighbours of $x$ map to $t_1$ and $t_2$.
Let $y$ be the out-neighbour that maps to $t_1$.
The vertex $y$ has another in-neighbour, $z$, which by injectivity maps to $t_0$.
Now, following the same argument as in the previous paragraph 
we have that, starting from $x$, the images of consecutive vertices
are $t_1, t_1, t_0, t_2, t_1, \ldots t_1, t_1, t_0, t_2, t_1$.
Again, $n \equiv 0$ (mod 4).

It now follows that an iot-injective homomorphism exists if and only
if  $n \equiv 0$ (mod 4).
\hfill$\Box$

\begin{cor}
For $t \geq 1$, let $B_{6t} \in \mathcal{B}$ have $6t$ vertices.  In any iot-injective homomorphism $f$ of $B_{6t}$ to $C_3^r$ we have $f(v_i)=f(v_j)$, when $i\equiv j$ (mod 6).  
\label{CorollaryB_18}
\end{cor}
\noindent{\bf Proof.}
This follows from the argument in Lemma \ref{hmm}.
\hfill$\Box$

\begin{thm}
The problem of deciding whether an oriented graph has an iot-injective homomorphism to $C_3^r$ is NP-complete.
\label{thm:C3}
\end{thm}

\noindent \textbf{Proof.}  
The transformation is from 3-colouring of connected graphs  \cite{Garey}.  
Suppose such a graph $G$ is given.  
Construct a graph $G^\prime$ as follows.
For each $x \in V(G)$, regard the edges incident with $x$ as being in 1--1 correspondence with the integers $1, 2, \ldots, \mathrm{deg}(x)$ so that it is meaningful to talk about the $i^{\textup{th}}$ edge incident with $x$. 
Replace every vertex $x\in V(G)$ with a copy $R_x$ of $B_{6 \cdot \mathrm{deg}(x)}$ where,
without loss of generality, the vertices $x_{6i} \in V(R_x),\ 0 \leq i \leq \mathrm{deg}(x)-1$ have
in-degree 2.
Suppose $xy$ is the $i^{\text{th}}$ edge incident with  $x$ and the $j^{\text{th}}$ edge incident with $y$.   
Construct an oriented path $P_{xy}$ by adding a new vertex $t_{xy}$ and joining each of 
$x_{6(i-1)}\in V(R_x)$ and $y_{6(j-1)}\in V(R_y)$ to it by adding a directed path of length two (the midpoint of each such directed path is a new vertex).
The transformation can be carried out in polynomial time.  We claim that $G$ is $3$-colourable if and only if $G'$ has an iot-injective homomorphism to $C_3^r$.

Suppose $G$ has a $3$-colouring $f:V(G) \rightarrow \{c_1,c_2,c_3\}$.  
For each vertex $x$, map the vertices $x_{0},x_{6},\ldots x_{6 \cdot \mathrm{deg}(x)}$ of $R_x$ to $f(x)$.   
By Corollary \ref{CorollaryB_18}, this partial mapping extends to an iot-injective homomorphism of $R_x$ to $C_3^r$.  
We claim that this mapping of the oriented cycles $R_x$ extends to the oriented paths $P_{xy}$, where $xy \in E(G)$.  Since adjacent vertices in $G$ must receive different colours, this mapping of the copies of $B_{6 \cdot \mathrm{deg}(x)}$ assigns the vertices $v_{0},v_{6},\ldots v_{6 \cdot \mathrm{deg}(x)}$ of $R_x$ a different image than it assigns the corresponding vertices of $R_y$. Suppose, without loss of generality, that the vertices $v_{0},v_{6},\ldots, v_{6 \cdot \mathrm{deg}(x)}$ of $R_x$ are mapped to $c_1$ and the corresponding vertices of $R_y$ are mapped to $c_2$.   The in-neighbours of the vertices in $R_x$ are mapped to $c_1$ and $c_3$, while the neighbours of the corresponding vertices in $R_y$ are mapped to $c_2$ and $c_1$.  The vertex $t_{xy}$ can be mapped to $c_3$ and the assignment extended to an iot-injective homomorphism of $P_{xy}$ to $C_3^r$.  This proves the claim, and completes the proof of the implication.

On the other hand, suppose $G^\prime$ has an iot-injective homomorphism to $C_3^r$.  Fix such a mapping. 
Then, for each $v \in V(G)$, 
the vertices  $v_{0},v_{6},\ldots v_{6 \cdot \mathrm{deg}(v)} \in V(R_v)$ all have the same image; assign this to be the colour of vertex $v$ of $G$. 

We claim that vertices $x$ and $y$ that are adjacent in $G$ are assigned different colours.
Suppose not.
By symmetry of $C_3^r$, assume both are assigned $c_1$.
Suppose also that $xy$ is the $i$-th edge incident with $x$ and the $j$-th edge incident with $y$.
Then,  the vertices $x_{6(i-1)} \in V(R_x)$ and $y_{6(j-1)} \in V(R_y)$ both
map to $c_1$.
By construction,  $x_{6(i-1)}$ has two in-neighbours in $R_x$ and one out-neighbour on the directed path to $t_{xy}$,
and similarly for $y_{6(j-1)}$.
Since both $x_{6(i-1)}$ and $y_{6(j-1)}$ map to $c_1$, in each case their in-neighbours must map to $c_1$ and $c_3$.
By injectivity, in each case their out-neighbour on the directed path to $t_{xy}$ must map to $c_2$.
Therefore $t_{xy}$ has 2 in-neighbours that map to $c_2$, which violates injectivity.
This proves the claim, and completes the proof.
\hfill $\square$

\medskip


Finally, we consider iot-injective homomorphism to $T_3^r$.
The following lemma can be proved similarly to Lemma  \ref{ios-gadget}.
The proof of Lemma  \ref{ios-gadget} relies only on injectivity on
in-neighbourhoods or out-neighbourhoods, and never both
at the same vertex.

\begin{lem}
Let $F$ be the oriented graph in Figure \ref{figF}.  Then
For $x \in \{t_0, t_1, t_2\}$, there exists an iot-injective homomorphism of $F$ to $T_3^r$ that maps $u$ to $x$, and any such homomorphism also maps $v$ to $x$.
\label{iot-gadget} 
\end{lem} 

The proof of the following theorem is similar to that of Theorem \ref{iosT_3^r} and is omitted. For details, see \cite{russell}.

\begin{thm}
The problem of deciding whether an oriented graph has an iot-injective homomorphism to $T_3^r$ is
NP-complete.
\end{thm}
  
\section{Colourings}

Recall that a (proper) \emph{oriented $k$-colouring} of an oriented graph $G$ is a homomorphism to a tournament on $k$ vertices.  
We therefore make the following definitions:
\begin{enumerate}
\item A \emph{proper ios-injective oriented $k$-colouring} of an oriented graph $G$ is an ios-injective homomorphism to an irreflexive tournament on $k$ vertices. 
\item An \emph{improper ios-injective oriented $k$-colouring} of an oriented graph $G$ is an ios-injective homomorphism to a reflexive tournament on $k$ vertices.
\item An \emph{improper iot-injective oriented $k$-colouring} of an oriented graph $G$ is an iot-injective homomorphism to a reflexive tournament on $k$ vertices.
\end{enumerate}
A proper iot-injective oriented $k$-colouring of a graph $G$ would be an iot-injective homomorphism to
an irreflexive tournament on $k$ vertices. Since tournaments have no directed 2-cycles, these are the same as proper ios-injective oriented $k$-colourings.

For each fixed integer $k$ and each injective colouring problem defined above, we will determine the complexity of deciding whether a given oriented graph $G$ has an injective colouring with $k$ colours.
The approach to proving NP-completeness is similar to that  for oriented colourings that are injective on in-neighbourhoods \cite{mrs,cobusthesis}: prove that it is NP-complete to decide the existence of an injective homomorphism of the given type to the tournament $U_m, \ m \geq 4$, that consists of a directed three cycle dominated by every vertex of a transitive tournament of size $m-3$, and then obtain the desired result as a corollary.
We consider the three situations in turn after establishing a useful lemma.

\begin{lem}
Let $G$ be an oriented graph such that $U_m$ is a subgraph of $G$.  For $\mathcal{P} \in \{\mathrm{ios,\ iot}\}$, if $G$ has a $\mathcal{P}$-injective homomorphism to a tournament $T$ (respectively, reflexive tournament $T^r$), then $U_m$ (respectively, $U_m^r$) is a subgraph of $T$. 
\label{ontoU_m}
\end{lem} 

\noindent \textbf{Proof.} 
The tournament $U_m$ has the property that every two different vertices have a common in-neighbour or a common out-neighbour.  Hence no two of its vertices can be assigned the same image by a $\mathcal{P}$-injective homomorphism.  Consequently, the image of $G$ must contain $U_m$.  
\hfill $\square$

\subsection{Proper ios-injective colourings}

\begin{thm}
For each fixed $m \geq 4$, the problem of deciding if a given oriented graph has an ios-injective homomorphism to $U_m$ is NP-complete.
\end{thm}

\noindent \textbf{Proof.} 
We first show that 
the problem of deciding whether a given oriented graph $G$ has an ios-injective homomorphism to $U_4$ is NP-complete. 
The transformation is from the problem of deciding if a given cubic graph is 3-edge-colourable~\cite{holyer}.  
Let $G$ be a given cubic graph.  Construct an oriented graph $G^\prime$ by replacing each edge $xy$ of $G$ by an oriented path $P_{xy}$ with vertices $x, v_1, v_2, v_3, v_4, y$ and arcs $xv_1, v_1v_2,$ $v_2v_3, v_3v_4,$ $yv_4$.  The transformation can be accomplished in polynomial time. We claim that $G$ is $3$-edge-colourable if and only if there is an ios-injective homomorphism of $G^\prime$ to $U_4$.

Suppose that $G$ is $3$-edge colourable.  Then, for each vertex $x$ of $G$, each of the colours $1, 2$, and $3$ appears on an edge incident with $x$.  An ios-injective homomorphism of $G^\prime$ to $U_4$ is obtained by mapping all vertices of $G$ to the vertex of out-degree $3$ in $U_4$, assigning the colour of the edge $xy$ to the vertices $v_1$ and $v_4$ of $P_{xy}$, and extending this pre-colouring to the vertices $v_2$ and $v_3$ of $P_{xy}$.

Suppose that there is an ios-injective homomorphism of $G^\prime$ to $U_4$.  Every vertex of $G$ has out-degree $3$ in $G^\prime$, so an ios-injective homomorphism of $G^\prime$ to $U_4$ must map it to the unique vertex of out-degree $3$ in $U_4$.  Similarly, the vertices $v_1, v_2, v_3,$ and $v_4$ in each oriented path $P_{xy}$ have positive in-degree in $G^\prime$, so an ios-injective homomorphism of $G^\prime$ to $U_4$ must map each of them to a vertex of the directed 3-cycle.  In any such mapping, $v_1$ and $v_4$ map to the same vertex, and the three out-neighbours of each vertex of $G$ (in $G^\prime$) map to different vertices of the 3-cycle.  Assigning each edge $xy$ of $G$ the image of the vertex $v_1$ (and $v_4$) in $P_{xy}$ gives a 3-edge-colouring of $G$.  

NP-completeness of ios-injective homomorphism to $U_m$ follows from NP-completeness of ios-injective homomorphism to $U_4$.  Given an instance $G$ of ios-injective homomorphism to $U_4$, construct $G^\prime$ by adding the new vertices belonging to $V^\prime = \{x_i: x \in V(G), i = 1, 2, \ldots, (m-4)-\mathrm{d}^-(x)\}$ and the arcs $\{x_i x: x \in V(G), i = 1, 2, \ldots, (m-4)-\mathrm{d}^-(x)\}$. Since $m$ is a constant, the transformation can be accomplished in polynomial time.  Each vertex of $G$ in $G^\prime$ has in-degree $m-4$ and therefore cannot map to the $m-4$ vertices of $U_m$ with in-degree less than $m-4$.  An ios-injective homomorphism of $G$ to $U_4$ can be extended to an ios-injective homomorphism of $G^\prime$ to $U_m$.
\hfill $\square$

\begin{cor}
Let $k$ be a fixed positive integer.  If $k \leq 3$, the problem of deciding if a given oriented graph $G$ has a proper ios-injective oriented $k$-colouring is Polynomial.
If $k \geq 4$, the problem of deciding if a given oriented graph $G$ has a proper ios-injective oriented $k$-colouring is NP-complete.  
\end{cor}

\noindent \textbf{Proof.}
An oriented graph $G$ has a proper ios-injective oriented $k$-colouring if and only if $G \cup U_k$ has an ios-injective homomorphism to $U_k$.
\hfill $\square$

\subsection{Improper ios-injective colourings}

\begin{thm}
For each fixed $m \geq 4$, the problem of deciding if a given oriented graph has an ios-injective homomorphism to  $U_m^r$ is NP-complete.
\label{ThmIOSU_k^r}
\end{thm}

\noindent \textbf{Proof.} The transformation is from the problem of deciding whether there exists an ios-injective homomorphism of a given oriented graph $G$ to $C_3^r$, which is NP-complete by
Theorem  \ref{C_3^rNPc}.
Suppose the oriented  graph $G$ is given.
We may assume that 
$\Delta^+(G) \leq 2$ and $\Delta^-(G) \leq 2$, otherwise $G$ cannot have an ios-injective homomorphism to $C_3^r$. 

Construct $G^\prime$ from $G$ as follows.  For each $x \in V(G)$, if $x$ has in-degree at most one in $G$, add a set of $m-2$ new vertices and arcs joining each of them to $x$.  If $x$ has in-degree two in $G$, do the same using a set of $m-3$ new vertices.  The transformation can be accomplished in polynomial time. We claim that $G$ has an ios-injective homomorphism to $C_3^r$ if and only if $G^\prime$ has an ios-injective homomorphism to $U_m^r$. 

An ios-injective homomorphism of $G$ to $C_3^r$ can clearly be extended to an ios-injective homomorphism of $G^\prime$ to $U_m^r$.

Suppose $f$ is an ios-injective homomorphism of  $G^\prime$ to $U_m^r$.  Since each vertex $x \in V(G)$  has in-degree at least   $m-1$ in $G^\prime$ and every vertex of $U_m^r$ not belonging to the copy of $C_3^r$ has in-degree at most $m-3$, the vertex $x$ must map to a vertex of the directed 3-cycle in $U_m^r$.  The restriction of $f$ to $V(G)$ is the desired mapping. 
\hfill $\square$

\begin{cor}
Let $k$ be a fixed integer.  If $k \leq 2$, the problem of deciding if a given oriented graph $G$ has an improper ios-injective oriented $k$-colouring is Polynomial.
If $k \geq 3$, the problem of deciding if a given oriented graph $G$ has an improper ios-injective oriented $k$-colouring is NP-complete.
\label{ImproperIOS}
\end{cor}

\noindent \textbf{Proof.}
When $k=3$ the transformation is from the problem of deciding whether there exists an ios-injective homomorphism of a given oriented graph $G$ to  $C_3^r$, which is NP-complete
by Theorem \ref{C_3^rNPc}.  Since there is no ios-injective homomorphism of $D_d$ 
(from Lemma \ref{LemmaD_d}) to $T_3^r$, an oriented graph $G$ has an improper ios-injective oriented $3$-colouring if and only if $G \cup D_6$ has an ios-injective homomorphism to $C_3^r$.

When $k \geq 4$, the transformation is from the problem of deciding whether there exists an ios-injective homomorphism of a given oriented graph $G$ to $U_k^r$.  Given an oriented graph $G$, the transformed instance of improper ios-injective oriented $k$-colouring is the oriented graph $G \cup U_k$.   The claim that this $G^\prime$ has an improper ios-injective oriented $k$-colouring if and only if $G$ has an ios-injective homomorphism to $U_k^r$ follows from Lemma \ref{ontoU_m}.
\hfill $\square$

\newpage
\subsection{Improper iot-injective colourings}

\begin{thm}
For each fixed $m \geq 4$, the problem of deciding if a given oriented graph has an improper iot-injective homomorphism to $U_m^r$ is NP-complete.
\end{thm}

\noindent \textbf{Proof.} 
The transformation is from the problem of deciding whether there exists an iot-injective homomorphism of a given oriented graph $G$ to $C_3^r$, which is NP-complete by Theorem \ref{C_3^rNPc}.  Given an oriented graph $G$, the transformed instance $G^\prime$ is constructed by starting with $G$ and proceeding as follows.  For each vertex $x \in V(G)$, add a copy of $T_{m-3}$ and arcs from each of its vertices to $x$.  Then for every vertex $t$ of each copy of $T_{m-3}$ that was added, add three vertices, $t_a,\ t_b,\ t_c$, and arcs from $t$  to each of them.  The oriented graph $G$ has an iot-injective $C_3^r$-colouring if and only if $G'$ has an iot-injective $U_m^r$-colouring.  
\hfill $\square$

\medskip
The proof of the following is identical to that of Corollary \ref{ImproperIOS}, except for replacing ``ios'' by ``iot''.

\begin{cor}
Let $k$ be a fixed integer.
If $k \leq 2$, the problem of deciding whether an oriented graph has an improper iot-injective $k$-colouring is Polynomial.
If $k \geq 3$, the problem of deciding whether an oriented graph has an improper iot-injective $k$-colouring is NP-complete.
\label{cor:np}
\end{cor}


For a given oriented graph $G$, we denote by $\chi_{\mathit{ios}}(G), \chi_{\mathit{ios}}^r(G)$ and $\chi_{\mathit{iot}}^r(G)$, the smallest number of colours in a proper ios-injective oriented colouring, an improper ios-injective oriented colouring, and an improper iot-injective oriented colouring of $G$, respectively.  The superscript ``$r$'' is used to designate the improper colourings because the target graph being reflexive is what allows adjacent vertices to be assigned the same colour.  
A project for future research is to find tight bounds for these parameters.
The upper bounds should be exponential in the in-degree and out-degree -- consider the disjoint union of all tournaments on a fixed number of vertices.
Weak upper bounds can be obtained using the methods in \cite{russell,mrs,mrs3,cobusthesis}.
Tight bounds and efficient algorithms for trees can be obtained as in \cite{russell,mrs,mrs3}.

\end{document}